\documentclass[a4paper, 11pt]{article}
\usepackage[cp1251]{inputenc}
\usepackage{amsfonts,amsmath,amssymb,graphicx}

\begin{document}

\begin{center}
{\bf Rational solutions of the fourth and fifth Painlev\'{e} hierarchies}
\end{center}

\begin{center}
\it A.\,A.\,Grigor'ev (Minsk, Belarus)
\end{center}
We will consider four hierarchies of higher order analogues of the fourth (P4) and fifth (P5)
Painleve equations. The necessary and sufficient conditions for having rational solutions will be
presented. Also the algorithm for obtaining such solutions will be described.

The hierarchies obtained in \cite{NY-paper} generalize the property of having an extended affine
Weyl group of B{\"a}cklund transformations. The hierarchy of (P4) is represented by the system of
ordinary differential equations
\begin{equation}\label{Nhierarchy}
\left\{
\begin{array}{rl}
  &f_0'=f_0(f_1-f_2+f_3-f_4+...+f_{n-2}-f_{n-1}) + \alpha_0\\
  &f_1'=f_1(f_2-f_3+f_4-f_5+...+f_{n-1}-f_{0}) + \alpha_1\\
  &...\\
  &f_{n-1}'=f_{n-1}(f_0-f_1+f_2-f_3+...+f_{n-3}-f_{n-2}) + \alpha_{n-1},\\
\end{array}
\right.
\end{equation}
where $\alpha_i$ are such complex parameters, that $\sum_{i=0}^{n-1} \alpha_i=h \neq 0$, $f_i$ are
unknown functions and $n$ is an odd number defining different members of the hierarchy. Further
description of (\ref{Nhierarchy}) and its properties connected to this paper could be found in
\cite{NY-book}, \cite{GrigorievP4} and \cite{ClarksonFilipuk}.

Another hierarchy of analogues of the (P4) was considered by Shabat \cite{Shabat}, Veselov
\cite{ShabVes} and Adler \cite{Adler}. It was obtained as the result of periodic closing of the
dressing chain for linear Schr{\"o}dinger equation.
\begin{equation} \label{AShierarchy}
\left\{
\begin{array}{rl}
  &g_0' + g_1'=g_1^2-g_0^2 + \alpha_0\\
  &g_1' + g_2'=g_2^2-g_1^2 + \alpha_1\\
  &...\\
  &g_{n-1}'+g_0'=g_0^2-g_{n-1}^2 + \alpha_{n-1},\\
\end{array}
\right.
\end{equation}
Here $\alpha_i$ are also such complex parameters that $\sum_{i=0}^{n-1} \alpha_i=h \neq 0$ and $n$
is odd. With the scaling change of variables in (\ref{Nhierarchy}) as also in (\ref{AShierarchy})
one can vary the value of $h$. Considering this we will further use $h=1$ without the loss of
generality.

As analogues of (P5) we will again consider the hierarchies of Noumi and Yamada from
\cite{NY-paper} with ($\sum_{i=0}^{2n+1} \alpha_i=h \neq 0$) and arbitrary natural $n$ defining the
number of a member in the hierarchy. After change of variables it takes form
\begin{equation}\label{mPNsymAC}
\left\{ \begin{array}{l}
 z f_i'= z f_i \Phi_i - A_{Mod(i,2)} f_{0} + \alpha_i C_{Mod(i,2)}, \: i=\overline{0,2n+1} \\
 \sum \limits_{r=0}^n f_{2r}=C_0 , \:\:
    \sum \limits_{s=0}^n f_{2s+1}=C_1
\end{array} \right.
\end{equation}
where
\begin{equation}\label{Oboz}
\begin{array}{c}
    \Phi_i^{(n)}= \sum \limits_{1\leq r \leq s \leq n} f_{i+2r-1}f_{i+2s} - \sum \limits_{1\leq r \leq s \leq n} f_{i+2r}f_{i+2s+1}, \:\:
    \\ A_0^{(n)} =\sum \limits_{r=0}^n \alpha_{2r}, \:\:
    A_1^{(n)} =\sum \limits_{s=0}^n \alpha_{2s+1}.
    \end{array}
\end{equation}
 Some properties of (\ref{mPNsymAC}) connected to
this paper are presented in \cite{NY-paper}, \cite{GrigorievPolotsk} and \cite{GrigorievP5}. We
will explain the correspondance of (\ref{mPNsymAC}) to the systems of ordinary differential
equations obtained by Shabat \cite{Shabat}, Veselov \cite{ShabVes} and Adler \cite{Adler} with
($\sum_{i=0}^{2n+1} \alpha_i=h \neq 0$).
\begin{equation} \label{AShierarchy5}
\left\{
\begin{array}{rl}
  &g_0' + g_1'=g_1^2-g_0^2 + \alpha_0\\
  &g_1' + g_2'=g_2^2-g_1^2 + \alpha_1\\
  &...\\
  &g_{2n+1}'+g_0'=g_0^2-g_{2n+1}^2 + \alpha_{2n+1},\\
\end{array}
\right.
\end{equation}
Without the loss of generality we will again assume for (\ref{mPNsymAC}) and (\ref{AShierarchy5})
that $h=1$.

When considering these hierarchies with real-valued parameters $\alpha_i$, there exist a
B{\"a}cklund transformation, which leads to the system with $0\leq\alpha_i\leq1$ e.g.
\cite{GrigorievP4}. In this case the following sufficient conditions on $\alpha_i$ for existing of
rational solutions are known: for the system (\ref{Nhierarchy}) with odd $n$
\begin{equation}\label{NhierarchyRatParam}
\begin{array}{cc}
    \begin{array}{cccccccl}
      (\alpha_0, & \alpha_1, & \alpha_2, & \alpha_3, & \alpha_4, &..,& \alpha_{n} &) \\ \hline
      (1, & 0, & 0, & 0, & 0,&..., & 0 &), \:  \\
      (\frac{1}{3}, & \frac{1}{3}, & \frac{1}{3}, & 0, & 0, &...,& 0 &), \: \\
      (..., & ..., & ..., & ..., & ...,&..., & ... &) \\
      (\frac{1}{n}, & \frac{1}{n}, & \frac{1}{n}, & \frac{1}{n}, & \frac{1}{n},&..., & \frac{1}{n} &). \:
    \end{array} &
    \begin{array}{cccccccl}
      (f_0, & f_1, & f_2, & f_3, & f_4, &..,& f_{n} &) \\ \hline
      (t, & 0, & 0, & 0, & 0,&..., & 0 &), \:  \\
      (\frac{t}{3}, & \frac{t}{3}, & \frac{t}{3}, & 0, & 0, &...,& 0 &), \: \\
      (..., & ..., & ..., & ..., & ...,&..., & ... &) \\
      (\frac{t}{n}, & \frac{t}{n}, & \frac{t}{n}, & \frac{t}{n}, & \frac{t}{n},&..., & \frac{t}{n} &). \:
    \end{array}
\end{array}
\end{equation}

For the system (\ref{mPNsymAC}) with $C_0 \neq 0$ and $C_1 \neq 0$ we also have
\begin{equation}\label{mPNsymACRatParam}
    \begin{array}{cccccccl}
      (\alpha_0, & \alpha_1, & \alpha_2, & \alpha_3, & \alpha_4, &..,& \alpha_{2n+1} &) \\ \hline
      (\alpha_0, & 1-\alpha_0, & 0, & 0, & 0,&..., & 0 &), \: \alpha_0 \in [0,1] \\
      (\alpha_0, & \frac{1}{2}(1-2 \alpha_0), & \alpha_0, & \frac{1}{2}(1-2 \alpha_0), & 0, &...,& 0 &), \: \alpha_0 \in [0,\frac{1}{2}] \\
      (..., & ..., & ..., & ..., & ...,&..., & ... &) \\
      (\alpha_0, & \frac{1}{n+1}-\alpha_0, & \alpha_0, & \frac{1}{n+1}-\alpha_0, & \alpha_0,&..., & \frac{1}{n+1}-\alpha_0 &), \: \alpha_0 \in [0,\frac{1}{n+1}]
    \end{array}
\end{equation}

And if in the rows of the tables (\ref{NhierarchyRatParam}) and (\ref{mPNsymACRatParam}) we admit
all possible permutations of all couples of zeroes this will describe the complete set of the
parameters $\alpha_i$, with which the hierarchies (\ref{Nhierarchy}) and (\ref{mPNsymAC}) admit
rational solutions. Considering the relations between the both P4-hierarchies
(\ref{Nhierarchy}),(\ref{AShierarchy}) and also between the both P5-hierarchies (\ref{mPNsymAC}),
(\ref{AShierarchy5}) the tables could be used for obtaining the conditions of having rational
solutions for the systems (\ref{AShierarchy}) and (\ref{AShierarchy5}).

Further we well consider two more hierarchies of the (P4). We well show that both of them admit the
solution in the form of $1/x$. These hierarchies are the one that was introduced by Kudryashov in
\cite{Kudr-hier} as the compatibility condition
\begin{equation}\label{Kudr-Comp-Cond}
    4 P A_x+2 A P_x+\lambda A y_x+ 2 \lambda y A_x - y y_x A+A_{xxx}+\lambda^2+2y_x A_x-A
    y_{xx}-\lambda^2 A_x-y^2A_x-\lambda y=0
\end{equation}
for the system of linear partial differential equations
\begin{equation}\label{Kudr-Pair}
\left\{
\begin{array}{c}
  \Psi_{xx}=(\lambda - y(x)) \Psi_x - P(x) \Psi \\
  \lambda \Psi_{\lambda} = A(x,\lambda) \Psi_x + B(x,\lambda)\Psi,
\end{array}
\right.
\end{equation}
where $A(x,\lambda)$ is a polynomial of $\lambda$. The degree of $A(x,\lambda)$ as a polynomial of
$\lambda$ corresponds to the member's number in the hierarchy. For example, the second member
($n=2$) has the form:
\begin{equation}\label{Kudr4}
    \begin{array}{c}
      (y''-2 x y-2 y^3-\beta)y^2 y^{(4)} - \frac{1}{2} y^2 \left(y^{(3)}\right)^2 + (2 y^2 + 8 y^3 y' + 4x y y'-y' y''+\beta y')y y^{(3)}-\\
       -\frac{4}{3} y (y'')^3 + \left( 3 x y^2 +3 \beta y - \frac{3}{2}y^4 + \frac{3}{2} (y')^2\right)(y'')^2 + (\beta y^4 - 2y'y^2-12(y')^2 y^3-2 \beta^2 y +\\
       +10 x y^5 - 3 \beta (y')^2 +10 y^7 - 4x y (y')^2 -4 \beta x y^2)y''+ 2(\beta-4 y^3) y^2 y' + ( 4 \beta x y +8 x y^4 + \\
       + \frac{3}{2} \beta^2 +12 \beta y^3)(y')^2 -\frac{10}{3} y^{10} -8x y^8 - 2\beta y^7 - 6x^2 y^6 - 2 \beta x y^5+\\
+ \left(\frac{1}{2} \beta^2 - 2 +9\delta - \frac{4}{3} x^3 \right) y^4 + \beta x y^2 + \frac{1}{3}
\beta^3 y =0,
    \end{array}
\end{equation}
In this equation $\beta$ and $\delta$ are complex parameters.

And the last considered here analogues of the (P4) were obtained by Gordoa, Joshi and Pickering in
\cite{GJP-hier} from nonisospectral scattering problem
\begin{equation}\label{GJP4}
    \left\{
    \begin{array}{l}
     L_{n,x}-2 K_n-(u+2 \frac{g_{n+1}}{g_n})L_n=g_n-2 \alpha_n \\
     L_n K_{n,x}+v L^{2}_{n}+K^{2}_{n} -L_{n,x}K_{n}(u+2 \frac{g_{n+1}}{g_n})L_n K_n=(\frac{1}{2}g_n - \alpha_n)^2-\frac{1}{4}\beta^2_n
    \end{array} \right. ,
\end{equation}
where
\begin{equation}
\begin{array}{lll}
   \left( \begin{array}{l} K_n \\ L_n \end{array} \right) & = &
   B^{-1} \left[ R^{n-1} U_x + g_{n}  \left( \begin{array}{l} 1 \\ 0 \end{array} \right) + \sum \limits_{i=0}^{n-2} (\frac{-g_{n+1}}{g_n})^{n-i-1} R^{i} U_x \right] +2(\frac{-g_{n+1}}{g_n})^{n-i-1} \left( \begin{array}{l} 1 \\ 0 \end{array} \right)
\end{array}
\end{equation}
\begin{equation} \label{R-oper}
\begin{array}{lll}
U= \left( \begin{array}{l} u \\ v \end{array} \right)=  \left( \begin{array}{l} u(x) \\ v(x)
\end{array} \right), & B = \left(
      \begin{array}{cc}
        0 & \partial_x \\
        \partial_x & 0 \\
      \end{array}
    \right) ,
& R=\frac{1}{2} \left(
                \begin{array}{cc}
                  \partial_x u \partial_x^{-1} & 2 \\
                  2 v + v_x \partial_x^{-1} & u+ \partial_x \\
                \end{array}
              \right).
\end{array}
\end{equation}
In this problem $\alpha_n$, $\beta_n$ and $g_i, \; i=\overline{0,n+1}$ are complex parameters among
which $g_{n-1}=0$ and $g_n \neq 0$. The number $n$ denotes the member's place in the resulting
hierarchy.

The last result shows that all the considered P4-hierarchies have the common property: they assume
the solution $1/x$ as also the original P4 does.

\end{document}